\documentclass[12pt]{amsart}

\usepackage{amsmath, amscd, amssymb, latexsym}
\usepackage{graphicx}
\usepackage{hyperref}
\usepackage{cleveref}
\usepackage{mathptmx}
\usepackage{esint}

\theoremstyle{plain}

\newtheorem{theorem}{Theorem}[section]
\newtheorem{lemma}[theorem]{Lemma}

\numberwithin{equation}{section}
\allowdisplaybreaks[3]

\crefname{theorem}{Theorem}{Theorems}
\crefname{lemma}{Lemma}{Lemmas}
\crefname{proposition}{Proposition}{Propositions}
\crefname{corollary}{Corollary}{Corollaries}
\crefname{question}{Question}{Questions}
\crefname{conjecture}{Conjecture}{Conjectures}
\crefname{definition}{Definition}{Definitions}
\crefname{remark}{Remark}{Remarks}
\crefname{section}{Section}{Sections}
\crefname{figure}{Figure}{Figures}
\crefname{theoremx}{Theorem}{Theorems}

\begin{document}

\title[]
{Overdetermined boundary problems for maximal graphs}

\author{Mijia Lai}
\address{School of Mathematical Sciences, Shanghai Jiao Tong University, Shanghai 200240, China}
\email{laimijia@sjtu.edu.cn}
\dedicatory{To the memory of Guanghao Hong
}

\maketitle

\begin{abstract}
We study overdetermined boundary problems for the maximal surface
equation with zero Dirichlet data and constant nonzero Neumann data
. For the exterior of a bounded domain with connected
boundary in $\mathbb{R}^n$, $n\geq2$, we prove that such domain must be 
the complement of a ball and the solution is radial. In dimension two,
we also show that a positive solution on a domain bounded by a
properly embedded noncompact connected curve forces the domain
to be a half-plane and the solution to be affine. The proofs
combine the divergence structure of the equation with the
induced geometry of maximal graphs. 
\end{abstract}

\date{}
\section{Introduction}

Let $\mathbb{L}^{n+1}$ denote the Lorentz--Minkowski space endowed
with the metric
\[
\langle \cdot,\cdot\rangle_{\mathbb L}
=
dx_1^2+\cdots+dx_n^2-dx_{n+1}^2.
\]
The graph of a function $u:\Omega\subset\mathbb R^n\to\mathbb R$
is spacelike if $|Du|<1$, and its mean curvature with respect to the Lorentz metric vanishes precisely
when $u$ satisfies the maximal surface equation
\begin{equation}\label{eq:intro-maximal}
\operatorname{div}
\left(
\frac{Du}{\sqrt{1-|Du|^2}}
\right)=0
\qquad\text{in }\Omega.
\end{equation}
This equation is the Euler--Lagrange equation of the Lorentzian
area functional \[
    \mathcal{A}(u)
    =\int_{\Omega}\sqrt{1-|Du|^2}\,dx,
\]
under compactly supported variations of the graph.
The strict spacelikeness condition $|Du|<1$ ensures that the induced
metric $g_{ij}=\delta_{ij}-u_i u_j$ is Riemannian and that the equation
is elliptic. Unlike minimal graphs in Euclidean space, spacelike solutions of
\eqref{eq:intro-maximal} locally maximize the Lorentzian area among spacelike
graph variations with fixed boundary values, owing to the strict
concavity of the area integrand $\sqrt{1-|Du|^2}$.
This variational property explains the terminology \emph{maximal}.

The classical Bernstein theorem of Calabi
and Cheng--Yau asserts that every smooth entire spacelike
maximal graph must be a hyperplane
\cite{Calabi1970,ChengYau1976}.
For domains with boundary, the geometry of the boundary
and the prescribed boundary data allow a much richer
family of solutions; see, for example,
Bartnik--Simon \cite{BartnikSimon1982}.

On a bounded region $K$, the functional
\[
\mathcal J(u)=\int_K (1-\sqrt{1-|Du|^2})\,dx
\]
differs from the negative Lorentzian area by a constant,
and its Euler--Lagrange equation is
\eqref{eq:intro-maximal}.

We orient $\partial K$ by the unit ourward normal $\nu$. To describe domain variations without assuming that
$\mathcal J(u)$ is finite, fix a large ball $B_R$ and
consider deformations supported strictly inside $B_R$.
Let
\[
K_t=(\operatorname{Id}+tX)(K),
\qquad
\mathcal J_R(t)
=
\int_{B_R\setminus\overline{K_t}}
\left(1-\sqrt{1-|Du_t|^2}\right)\,dx,
\]
where $u_t$ solves the maximal surface equation, vanishes
on the moving boundary, and has fixed Dirichlet data on
the outer boundary.
For a smooth family of spacelike solutions, the shape
derivative is
\[
\mathcal J_R'(0)
=
\int_{\partial K\cap B_R}
\left(
\frac{1}{\sqrt{1-|Du|^2}}-1
\right)
X\cdot\nu\,d\sigma.
\]
The first variation of the enclosed volume is
\[
    \left.\frac{d}{dt}\right|_{t=0}|K_t|
    =\int_{\partial K}X\cdot\nu\,d\sigma.
\]
Thus, critical point under volume-preserving deformations implies there exists a Lagrange multiplier $\lambda$ such that
\[
    \frac{1}{\sqrt{1-|Du|^2}}-1=\lambda
    \qquad\text{on }\partial K.
\]

Hence we are led to the overdetermined problem
\begin{equation}\label{eq:intro-overdetermined}
\begin{cases}
\displaystyle
\operatorname{div}
\left(
\frac{Du}{\sqrt{1-|Du|^2}}
\right)=0
&\text{in }\Omega,\\[2mm]
u=0,\qquad \partial_\nu u=\kappa
&\text{on }\partial\Omega,
\end{cases}
\qquad 0<\kappa<1,
\end{equation}
where $\Omega=\mathbb{R}^n\setminus K$ and  $\nu$ denotes the unit normal
pointing into $\Omega$.
Geometrically, the graph has boundary in the horizontal
hyperplane and meets that hyperplane at a constant
hyperbolic angle.

The study of overdetermined elliptic problems originates
in Serrin's symmetry theorem \cite{Serrin1971}, which
establishes spherical symmetry for positive solutions
of a broad class of elliptic equations on bounded
domains with constant Dirichlet and Neumann data.
For exterior domains, Reichel
\cite{Reichel1996,Reichel1997} developed a version of
the moving plane method adapted to the unbounded
setting. Sirakov \cite{Sirakov2001} subsequently
extended this theory, allowing, in particular,
different constant boundary data on different
components of a bounded obstacle.
These results show that suitable overdetermined data,
together with appropriate conditions at infinity,
force the obstacle to be a ball.

For the corresponding Euclidean mean curvature equation,
the planar problem admits a particularly strong
classification. Cui \cite{Cui2022} proved that a positive
minimal graph with connected boundary, zero Dirichlet
data, and constant normal derivative must be a
half-plane or a portion of a vertical catenoid.
Equivalently, its domain is a half-plane or the
exterior of a disk. The proof uses analytic reflection
and classification results for complete minimal
surfaces.

Our first main result investigates (\ref{eq:intro-overdetermined}) on an exterior domain of a bounded set. 

\begin{theorem}
\label{thm:exterior-rigidity}
Let $n\ge2$ and let $K\subset\mathbb R^n$ be a bounded
domain with $C^{3,\alpha}$ connected boundary,
where $0<\alpha<1$.
Suppose that
$
\Omega=\mathbb R^n\setminus\overline K
$
is connected.
Let $u\in C^3(\Omega)\cap C^2(\overline\Omega)$
satisfy the overdetermined problem (\ref{eq:intro-overdetermined}).
Then there exist $x_0\in\mathbb R^n$ and $R>0$ such that
$
K=B_R(x_0)
$
and
\[
u(x)
=
\int_R^{|x-x_0|}
\frac{D}{\sqrt{r^{2n-2}+D^2}}\,dr,
\qquad
D=\frac{\kappa R^{n-1}}{\sqrt{1-\kappa^2}}.
\]
\end{theorem}

An important analytic input for this theorem is the asymptotic expansion of solutions to the maximal surface equation on exterior domains obtained by Hong and
Yuan \cite{HongYuan2021}. Let $\Omega$ be as above, and let $u$
be a smooth strictly spacelike solution of the maximal surface equation in $\Omega$. The divergence structure implies that the quantity \[
    Q[u]
    :=
    \int_{\partial B_R}
    \frac{Du\cdot\nu_R}{\sqrt{1-|Du|^2}}\,d\sigma
\]
is independent of $R$ whenever
$\overline{K}\subset B_R$. Following the terminology of Hong and Yuan,
the normalization quantity
$\operatorname{Res}[u]
    :=
    \frac{Q[u]}{\gamma_n}$
is called the \emph{residue} of $u$, where \[ \gamma_n
    :=
    \begin{cases}
        2\pi, & n=2\\
        (n-2)|\mathbb{S}^{n-1}|, & n\geq 3.
    \end{cases} 
    \]
Set
\[
\rho_a(x):=\sqrt{|x|^2-(a\cdot x)^2}.
\]
the asymptotic expansion of Hong--Yuan asserts that there exists $a\in \mathbb{R}^n, |a|<1$ and $b\in \mathbb{R}$ such that
\begin{equation}\label{eq: expansion}
u(x)=
\begin{cases}
a\cdot x+b
+(1-|a|^2)\operatorname{Res}[u]\log\rho_a(x)
+O_k\!\left(|x|^{-1}\log|x|\right),
& n=2,\\[2mm]
a\cdot x+b
-(1-|a|^2)\operatorname{Res}[u]\rho_a(x)^{2-n}
+O_k\!\left(|x|^{1-n}\right),
& n\geq 3,
\end{cases}
\end{equation}
for every integer $k\geq 0$.
Here the notation
$
R(x)=O_k\!\left(|x|^\mu(\log|x|)^\ell\right)
$
 means that, for every multi-index $\alpha$ with $|\alpha|\leq k$,
\[
|D^\alpha R(x)|
\leq C_\alpha |x|^{\mu-|\alpha|}
(\log|x|)^\ell
\]
for all sufficiently large $|x|$.

In dimension two, we have absorbed the explicit correction
term of order $|x|^{-1}\log|x|$ in Hong--Yuan's expansion
into the remainder.
When $a=0$, their sharper expansion gives
\[
u(x)=b+\operatorname{Res}[u]\log|x|
+O_k(|x|^{-1}).
\]

We next turn to domains with noncompact boundary. The geometry of overdetermined problems in this setting
is substantially richer. A central question was formulated by
Berestycki, Caffarelli, and Nirenberg \cite{BCN1997}, who conjectured
that a smooth domain $\Omega\subset\mathbb{R}^n$ with connected
complement admitting a bounded positive solution of
\[
    \Delta v+f(v)=0 \quad \text{in }\Omega,
    \qquad
    v=0,\quad \partial_\nu v=\mathrm{constant}\neq0
    \quad \text{on }\partial\Omega,
\]
for a Lipschitz nonlinearity $f$, must be a ball, a half-space,
a generalized cylinder, or the complement of one of these domains. The BCN conjecture is false in general.
Sicbaldi \cite{Sicbaldi2010} constructed counterexamples in dimensions
$n\geq3$ by perturbing straight cylinders periodically.
Ros, Ruiz, and Sicbaldi \cite{RosRuizSicbaldi2020} subsequently
constructed nonradial exterior domains supporting bounded positive
solutions for suitable nonlinearities in every dimension $n\geq2$,
thereby also disproving the conjecture in the plane.

For complements of unbounded obstacles, nontrivial examples exist
even in the harmonic case $f\equiv0$.
Fall, Minlend, and Weth \cite{FallMinlendWeth} constructed exceptional
domains in $\mathbb{R}^3$ as complements of nontrivial periodic
perturbations of a straight cylinder.
Minlend, Weth, and Wu \cite{MinlendWethWu} obtained analogous
constructions in every dimension $n\geq4$.
These domains have connected noncompact boundary and support positive
harmonic functions with zero Dirichlet data and constant nonzero
Neumann data. The corresponding functions have logarithmic growth
in the three-dimensional construction and are bounded in dimensions
$n\geq4$. Thus, in higher dimensions, connectedness of the boundary
does not force even a harmonic overdetermined domain to be a
half-space or the complement of a straight circular cylinder.
These examples also suggest that a similarly unrestricted
classification for maximal graphs should not be expected, although
they do not themselves provide counterexamples for the maximal
surface equation.

In dimension two, however, strong rigidity results remain available
when the boundary is connected and noncompact.
Ros, Ruiz, and Sicbaldi \cite{RosRuizSicbaldi2017} proved that a planar
domain with unbounded connected boundary supporting a bounded
positive solution of the semilinear overdetermined problem above
must be a half-plane, and the solution is one-dimensional.

Several rigidity results are also available
for geometric quasilinear equations in dimension two. In addition to the
aforementioned result of Cui \cite{Cui2022}, a recent contribution
is due to Lian and Sicbaldi \cite{LianSicbaldi2025}, who studied
the capillary overdetermined problem
\begin{align} \notag
\begin{cases}
    \operatorname{div}
    \left(\frac{Dv}{\sqrt{1+|Dv|^2}}\right)+f(v)=0
    \quad  \text{in }\Omega, \\
     v>0 \quad \text{in }\Omega,
    \qquad
    v=0,\quad \partial_\nu v=\kappa
    \quad \text{on }\partial\Omega.
\end{cases}
\end{align}
For planar domains with unbounded connected boundary, they obtained
half-plane rigidity and one-dimensionality of the solution under
boundedness of both $v$ and $Dv$, provided that $f$ admits a
nonpositive primitive $F$ satisfying
$
    F(0)\geq \frac{1}{\sqrt{1+\kappa^2}}-1.
$

Our second result establishes a rigidity statement in the setting of noncompact domains with connected unbounded boundary in dimension two.

\begin{theorem}
\label{thm:half-plane-rigidity}
Let $\Gamma\subset\mathbb{R}^2$ be a properly embedded
$C^{3,\alpha}$ curve diffeomorphic to $\mathbb{R}$, where
$\alpha\in(0,1)$, and let $\Omega$ be one of the two connected
components of $\mathbb{R}^2\setminus\Gamma$. Let $\kappa\in(0,1)$ be a constant.
Suppose that 
$u\in C^3(\Omega)\cap C^2(\overline{\Omega})$ satisfies the overdetermined problem \eqref{eq:intro-overdetermined}, 
then, after a rigid motion,
$\Omega=\{(x_1,x_2)\in\mathbb{R}^2:x_2>0\}$, and $u(x_1,x_2)=\kappa x_2.$
\end{theorem}

The proofs of both theorems combine the divergence structure of
the equation with the geometry of maximal graphs, making essential
use of the induced Riemannian metric. For Theorem~\ref{thm:exterior-rigidity}, integral identities and geometric inequalities
lead to rigidity through their equality cases. This strategy
is reminiscent of Weinberger's approach to Serrin's problem
\cite{Weinberger1971} and of the integral methods underlying
Alexandrov-type rigidity results \cite{Reilly1977,Ros1987}. For Theorem~\ref{thm:half-plane-rigidity}, a key step is
to establish the uniform spacelikeness ($|Du|\leq 1-\delta, \delta>0$), which allows us to exploit
global conformal coordinates and reduce the rigidity argument to the analysis of harmonic functions.

We present the proof of Theorem~\ref{thm:exterior-rigidity} in Section 2 and the proof of Theorem~\ref{thm:half-plane-rigidity} in Section 3.

\section{Rigidity for exterior domains of bounded obstacles}
\label{sec:exterior-rigidity}

Throughout this section,
 we write
\[
\Omega=\mathbb R^n\setminus\overline K,
\qquad
\Gamma=\partial K,
\qquad
W=\sqrt{1-|Du|^2},
\]
and orient $\Gamma$ by the unit normal $\nu$ pointing
from $K$ into $\Omega$. We denote the Euclidean area of $\Gamma$ by
$A=|\Gamma|$, the volume of $K$ by $V=|K|$. The mean curvature $H$ is the sum of the principal
curvatures with respect to $\nu$, so that
$H=(n-1)/R$ on $\partial B_R$. $\nu_R$ denotes the unit outward normal of $B_R$.

In view of the maximal surface equation \eqref{eq:intro-overdetermined}, the divergence theorem in $B_R\setminus\overline K$ gives
\begin{equation}\label{eq:horizontal-end-scalar-flux}
     \gamma_n\operatorname{Res}[u]=\int_{\partial B_R}\frac{Du\cdot\nu_R}{W}\,d\sigma
    =\int_{\partial K} \frac{u_{\nu}}{\sqrt{1-|Du|^2}} d\sigma= \frac{\kappa}{\sqrt{1-\kappa^2}} A.
\end{equation}
For simplicity, we set
\[
s_0:=\frac{\kappa}{\sqrt{1-\kappa^2}}.
\]

The first step is to show that the overdetermined
boundary conditions eliminate the linear part allowed
by the general exterior asymptotic expansion \eqref{eq: expansion}.
\begin{lemma}
\label{lem:horizontal-end}
Under the assumptions of Theorem~\ref{thm:exterior-rigidity},
the asymptotic slope $a$ of $u$ in the expansion (\ref{eq: expansion}) vanishes. 

If $n=2$, there exists $b\in\mathbb{R}$ such that
\begin{equation}\label{eq:horizontal-asymptotics-two}
    u(x)
    =
    \operatorname{Res}[u]\log|x|
    +b+O_k(|x|^{-1}).
\end{equation}
In this case, $u>0$ in $\Omega$.

If $n\geq3$, there exists $c>0$ such that
\begin{equation}\label{eq:horizontal-asymptotics-higher}
    u(x)
    =
    c-\operatorname{Res}[u]|x|^{2-n}
    +O_k(|x|^{1-n}).
\end{equation}
In this case, $0<u<c$ in $\Omega$.
The expansions hold for every nonnegative integer $k$.
\end{lemma}

\begin{proof}
Define
\[
    T_{ij}=W\delta_{ij}+\frac{u_i u_j}{W}.
\]
Each row of $T$ has zero divergence, since
\[
    \sum_{j=1}^n\partial_jT_{ij}
    =
    \partial_iW
    +\frac1W\sum_{j=1}^n u_{ij}u_j
    +u_i\operatorname{div}\left(\frac{Du}{W}\right)
    =0.
\]
On $\partial K$,
\[
    T\nu
    =
    \left(\sqrt{1-\kappa^2}+\frac{\kappa^2}{\sqrt{1-\kappa^2}}\right)\nu
    =
    \frac{\nu}{\sqrt{1-\kappa^2}}.
\]
The divergence theorem yields
\begin{equation}\label{eq:horizontal-end-vector-flux}
    \int_{\partial B_R}T\nu_R\,d\sigma
    =
    \frac1{\sqrt{1-\kappa^2}}\int_{\partial K}\nu\,d\sigma
    =0.
\end{equation}

By Hong and Yuan's asymptotic expansion \eqref{eq: expansion}, there exists
$a\in\mathbb{R}^n$, $|a|<1$, such that
\[
    Du=a+h,\qquad h=O(|x|^{1-n}).
\]
Write
\[
    W_a=\sqrt{1-|a|^2},
    \qquad
    v=u-a\cdot x,
    \qquad
    B=\frac{Du}{W}.
\]
Thus $h=Dv$, and
\[
    T-a\otimes B
    =
    WI+\frac{h\otimes(a+h)}{W}.
\]
Using
\[
    W=W_a-\frac{a\cdot h}{W_a}+O(|h|^2),
    \qquad
    \frac1W=\frac1{W_a}+O(|h|),
\]
we obtain
\begin{equation}\label{eq:horizontal-end-expansion}
    T-a\otimes B
    =
    W_aI
    +\frac{h\otimes a-(a\cdot h)I}{W_a}
    +O(|h|^2).
\end{equation}

Choose a smooth function $\widetilde v$ on
$\mathbb{R}^n$ which agrees with $v$ outside a fixed large ball,
using a cutoff near $\overline K$.
For each $i$, the vector field
$
    F^{(i)}_j
    =
    a_j\partial_i\widetilde v
    -\delta_{ij}\sum_{k=1}^n a_k\partial_k\widetilde v
$ is divergence free. Since $D\widetilde v=h$ on $\partial B_R$ for large $R$, it follows that each row of $h\otimes a-(a\cdot h)I$ is also divergence free, and consequently
\begin{equation}\label{eq:horizontal-end-linear-flux}
    \int_{\partial B_R}
    \bigl[h\otimes a-(a\cdot h)I\bigr]\nu_R\,d\sigma
    =0.
\end{equation}

Integrating \eqref{eq:horizontal-end-expansion} against $\nu_R$
and using \eqref{eq:horizontal-end-scalar-flux},
\eqref{eq:horizontal-end-vector-flux} and \eqref{eq:horizontal-end-linear-flux}, we find
\[
    -\gamma_n\operatorname{Res}[u]\,a
    =
    \int_{\partial B_R}O(|h|^2)\nu_R\,d\sigma.
\]
Consequently,
\[
    \gamma_n\operatorname{Res}[u]\,|a|
    \leq
    C\int_{\partial B_R}|h|^2\,d\sigma
    \leq C R^{1-n}.
\]
Letting $R\to\infty$ and using
$\operatorname{Res}[u]>0$, we conclude that $a=0$.

The expansions
\eqref{eq:horizontal-asymptotics-two} and
\eqref{eq:horizontal-asymptotics-higher} now follow from
\eqref{eq: expansion}.

For $n=2$, positivity of the residue implies
$u(x)\to+\infty$ as $|x|\to\infty$.
The maximum principle on bounded truncations, followed by the
strong maximum principle, gives $u>0$ in $\Omega$.

For $n\geq3$, suppose that $c\leq0$.
Since $u=0$ on $\partial K$ and $u(x)\to c$ at infinity,
the maximum principle on bounded truncations would give
$u\leq0$ in $\Omega$.
This contradicts $\partial_\nu u=\kappa>0$ on $\partial K$.
Hence $c>0$, and the maximum principle gives
$0<u<c$ throughout $\Omega$.
\end{proof}

\medskip

Next, we shall utilize the geometric property of the maximal surface equation. We equip $\Omega$ with the metric induced by the maximal graph,
$
    g_{ij}=\delta_{ij}-u_i u_j.
$
We will reserve $D$ for Euclidean derivative, and all notions with a subscript $g$ such as $\nabla_g$ are taken with
respect to $g$.

We list several geometric identities needed below.
\[
    g^{ij}=\delta_{ij}+\frac{u_i u_j}{W^2},
    \qquad
    \det g=W^2
    \qquad \nabla_g u=\frac{Du}{W^2}.
\]
Setting $
    s:=|\nabla_g u|_g=\frac{|Du|}{W}$, therefore 
    \[
    W=\sqrt{1-|Du|^2}= \frac1{\sqrt{1+s^2}}.   \]
From the maximal surface equation, it is easy to see that $u$ is harmonic with respect to $g$:
\begin{equation}\label{eq:height-g-harmonic}
    \Delta_g u
    =
    \frac1W\operatorname{div}\left(\frac{Du}{W}\right)
    =0.
\end{equation}

With the convention
$
    \mathrm{II}_{ij}=\frac{u_{ij}}{W}
$
for the second fundamental form in Lorentzian space, the Christoffel symbols and the Hessian of the height function are
\[
    \Gamma^k_{ij}=-\frac{u_k u_{ij}}{W^2},
    \qquad
    (\nabla_g^2u)_{ij}
    =\frac{u_{ij}}{W^2}
    =\frac{\mathrm{II}_{ij}}{W}.
\]
Since the graph is maximal, $\operatorname{tr}_g\mathrm{II}=0$.
The Gauss equation therefore gives
\begin{equation}\label{eq:maximal-ricci}
    \operatorname{Ric}_g(X,Y)
    =
    \sum_{j=1}^n
    \mathrm{II}(X,e_j)\mathrm{II}(Y,e_j),
\end{equation}
where $\{e_j\}_{j=1}^n$ is any $g$-orthonormal frame.

For later use, let $\Sigma$ be a regular level set of $u$, oriented
by the Euclidean unit normal
$
    \nu=\frac{Du}{|Du|}.
$
Its $g$-unit normal pointing toward increasing $u$ is
$
    \eta=\frac{\nabla_g u}{s}=\frac{\nu}{W}.
$
If $H=\operatorname{div}\nu$ denotes the sum of the Euclidean
principal curvatures of $\Sigma$, then its mean curvature in
$(\Omega,g)$, with the same trace convention, is
\[
    H_g=\operatorname{div}_g\eta=\frac{H}{W}.
\]
Consequently, \eqref{eq:height-g-harmonic} implies
\[
    0=\operatorname{div}_g(\nabla_g u)=\operatorname{div}_g(s\eta)
    =\eta(s)+sH_g.
\]
Therefore, we have
\begin{equation}\label{eq:level-normal-s}
    \eta(s)=-\frac{sH}{W}.
\end{equation}

\begin{lemma}\label{lem:gradient-transform} 
For $s>0$, define
\[
    F_n(s)=
    \begin{cases}
        \displaystyle
        \log\frac{s}{1+\sqrt{1+s^2}},
        & n=2,\\[3mm]
        \displaystyle
        \int_0^s
        \frac{\tau^{-1/(n-1)}}{\sqrt{1+\tau^2}}\,d\tau,
        & n\geq3.
    \end{cases}
\] 
Plugging $
    s:=|\nabla_g u|_g=\frac{|Du|}{W}$, then on $\{s>0\}$,
\begin{equation}\label{eq:gradient-transform-laplacian}
    \Delta_g F_n(s)
    =
    \frac{F_n'(s)}{s}
    \left(
        |\nabla_g^2u|_g^2
        -\frac{n}{n-1}|\nabla_g s|_g^2
    \right)
    \geq0.
\end{equation}
If $n=2$, equality holds in
\eqref{eq:gradient-transform-laplacian}.
If $n\geq3$, setting $F_n(0)=0$ gives a continuous weakly
subharmonic function on $\Omega$.
If $n=2$, for every $M>0$, the function
$\max\{F_2(s),-M\}$, assigned the value $-M$ at $s=0$,
is weakly subharmonic on $\Omega$.

Moreover, along every regular level set of $u$,
\begin{equation}\label{eq:normal-transform}
    \eta\bigl(F_n(s)\bigr)
    =-H\,s^{(n-2)/(n-1)}.
\end{equation}
\end{lemma}

\begin{proof}
We first work on $\{s>0\}$.
Since $\mathrm{II}=W\nabla_g^2u$ and
$
    \nabla_g^2u(\nabla_g u,\cdot)=s\,ds,
$
equation \eqref{eq:maximal-ricci} gives
\[
    \operatorname{Ric}_g(\nabla_g u,\nabla_g u)
    =
    W^2s^2|\nabla_g s|_g^2
    =
    \frac{s^2}{1+s^2}|\nabla_g s|_g^2.
\]
Applying the Bochner formula to the harmonic function $u$, we obtain
\begin{align*}
    \frac12\Delta_g s^2
    &=
    |\nabla_g^2u|_g^2
    +\operatorname{Ric}_g(\nabla_g u,\nabla_g u),\\
    s\Delta_g s
    &=
    |\nabla_g^2u|_g^2
    -\frac{|\nabla_g s|_g^2}{1+s^2}.
\end{align*}

For both definitions of $F_n$,
\[
    F_n'(s)
    =
    \frac{s^{-1/(n-1)}}{\sqrt{1+s^2}}>0,
    \qquad
    \frac{F_n''(s)}{F_n'(s)}
    =
    -\frac{1}{(n-1)s}-\frac{s}{1+s^2}.
\]
The chain rule thus yields
\begin{align*}
    \Delta_g F_n(s)
    &=
    F_n'(s)\Delta_g s
    +F_n''(s)|\nabla_g s|_g^2\\
    &=
    \frac{F_n'(s)}{s}
    \left(
        |\nabla_g^2u|_g^2
        -\frac{n}{n-1}|\nabla_g s|_g^2
    \right).
\end{align*}
By the refined Kato inequality, the term in the parentheses is nonnegative. When $n=2$, it vanishes. 

To handle critical points, suppose first that $n\geq3$.
The defining integral converges at zero, and $F_n(s)\to0$
as $s\to0$.
For every $\varepsilon>0$, the function
$\max\{F_n(s),\varepsilon\}$ is locally constant near
$\{s=0\}$ and weakly subharmonic on $\{s>0\}$.
It is therefore weakly subharmonic throughout $\Omega$.
Letting $\varepsilon\downarrow0$, with locally uniform convergence,
proves the assertion for $F_n(s)$.
For $n=2$, the same argument applies directly to
$\max\{F_2(s),-M\}$, since $F_2(s)\to-\infty$ as $s\to0$.

Finally, using \eqref{eq:level-normal-s} and
$F_n'(s)=W s^{-1/(n-1)}$, we find
\[
    \eta\bigl(F_n(s)\bigr)
    =
    F_n'(s)\eta(s)
    =
    W s^{-1/(n-1)}
    \left(-\frac{sH}{W}\right)
    =
    -H s^{(n-2)/(n-1)}.
\]
\end{proof}

\medskip

\begin{proof}[Proof of Theorem\ref{thm:exterior-rigidity}]

We divide the proof into two cases. 
\medskip

{\bf 1. The planar case}

Assume $n=2$.
By Lemma~\ref{lem:horizontal-end},
\[
u= \operatorname{Res}[u]\log|x|+b+O_1(|x|^{-1})
.
\]
It follows that
$
s=\frac{\operatorname{Res}[u]}{|x|}+O(|x|^{-2})
$
and hence
\[
F_2(s)
=
-\log|x|+\log\frac {\operatorname{Res}[u]}{2}+O(|x|^{-1}).
\]
Thus
$
G:=F_2(s)+\frac {u}{\operatorname{Res}[u]}
$
is bounded near infinity and equals $F_2(s_0)$
on $\Gamma$.

We claim that
\begin{equation}\label{eq:planar-comparison}
G\le F_2(s_0)
\qquad\text{in }\Omega.
\end{equation}
Indeed fix $\varepsilon>0$.
Since $u\to+\infty$ and $G$ is bounded near infinity,
$
G-\varepsilon u-F_2(s_0)<0
$
on every sufficiently large sphere.
On $\Gamma$, the same expression is zero.

To justify the comparison at critical points,
fix such a sphere $\partial B_R$ and choose $M$
large enough that the truncation
$\max\{F_2(s),-M\}$ agrees with $F_2(s)$ on both
boundary components of $\Omega\cap B_R$.
This is possible because $s>0$ on $\Gamma$ and
on sufficiently large spheres.
By Lemma~\ref{lem:gradient-transform},
\[
\max\{F_2(s),-M\}
+\left(\frac{1}{\operatorname{Res}[u]}-\varepsilon\right)u
-F_2(s_0)
\]
is weakly subharmonic.
The maximum principle gives the desired comparison
on $\Omega\cap B_R$.
Exhausting $\Omega$ and then letting
$\varepsilon\downarrow0$ proves
\eqref{eq:planar-comparison}.

It follows from \eqref{eq:normal-transform} that
\[
-H+\frac{s_0}{\operatorname{Res}[u]}\le0
\qquad\text{on }\Gamma.
\]
Since $Q=s_0A$ and $\operatorname{Res}[u]=Q/(2\pi)$,
\begin{align} \label{eq:planar-curvature-bound}
H\ge\frac{2\pi}{A}.
\end{align}
The connected boundary $\Gamma$ is a simple
closed curve, so the turning-angle formula gives
\[
\int_\Gamma H\,dS=2\pi.
\]
Equality must therefore hold everywhere in
\eqref{eq:planar-curvature-bound}.
Thus $\Gamma$ has constant positive curvature
and $K$ is a disk.

\medskip
{\bf 2. The higher-dimensional case}

Assume $n\ge3$ and abbreviate
\[
\beta=\frac1{n-1},
\qquad
F=F_n,
\qquad
F_0=F(s_0).
\]
By Lemma~\ref{lem:horizontal-end}, $0<u<c$ and
$u\to c$ at infinity.

The function
$
Z=F(s)-\frac{F_0}{c}(c-u)
$
is weakly subharmonic, vanishes on $\Gamma$,
and tends to zero at infinity.
The maximum principle on an exhaustion gives
\begin{equation}\label{eq:higher-comparison}
F(s)\le\frac{F_0}{c}(c-u)
\qquad\text{in }\Omega.
\end{equation}
Taking the inward derivative on $\Gamma$ yields
\[
-Hs_0^{1-\beta}+\frac{F_0s_0}{c}\le0.
\]
Consequently,
\begin{equation}\label{eq:mean-curvature-lower-bound}
H\ge h_0,
\qquad
h_0:=\frac{F_0s_0^\beta}{c}>0.
\end{equation}

The asymptotic expansion implies that
\[
J:=\int_\Omega(1-W)\,dx<\infty.
\]
Multiplying the maximal surface equation by $u$ and integrating
over $\Omega\cap B_R$, we find
\[
\int_{\Omega\cap B_R}\frac{|Du|^2}{W}\,dx
=
\int_{\partial B_R}
u\,\frac{Du\cdot \nu_R}{W}\,dS.
\]
Letting $R\to\infty$ gives
\begin{equation}\label{eq:height-flux-energy}
\int_\Omega\frac{|Du|^2}{W}\,dx=cQ.
\end{equation}

Consider the divergence-free tensor
\[
S=I-T
=(1-W)I-\frac{Du\otimes Du}{W}.
\]
Its trace is
\[
\operatorname{tr}S
=
n(1-W)-\frac{|Du|^2}{W}.
\]
On $\Gamma$,
\[
S\nu
=
-\left(\frac1{\sqrt{1-\kappa^2}}-1\right)\nu.
\]
Moreover, $S=O(|Du|^2)$ at infinity, so
\[
\int_{\partial B_R}x\cdot S\nu_R\,dS
=O(R^{2-n})\longrightarrow0.
\]
Applying the divergence theorem to $Sx$ and
using \eqref{eq:height-flux-energy}, we obtain
\begin{equation}\label{eq:pohozaev-exterior}
nJ-cQ
=
\left(\frac1{\sqrt{1-\kappa^2}}-1\right)
\int_\Gamma x\cdot\nu\,d\sigma
=
nV\left(\sqrt{1+s_0^2}-1\right).
\end{equation}

Define
\[
j(s)=\frac1{1+\sqrt{1+s^2}},
\qquad
I(s_0)=\int_0^{s_0}j(s)F'(s)\,ds.
\]
The function $j$ is decreasing.
For almost every $t\in(0,c)$, the level set
$\Sigma_t=\{u=t\}$ is smooth and compact.
The divergence theorem gives
\begin{equation}\label{eq:level-flux}
\int_{\Sigma_t}s\,d\sigma=Q.
\end{equation}

Since
\[
\frac{1-W}{|Du|}
=
\frac{s}{1+\sqrt{1+s^2}}
=sj(s),
\]
the Euclidean coarea formula yields
\[
J=\int_0^c\int_{\Sigma_t}sj(s)\,d\sigma\,dt.
\]
The critical set causes no additional contribution,
because $1-W=0$ wherever $Du=0$.

For $t\in[0,c]$, define $\bar s(t)\in[0,s_0]$ by
\[
F(\bar s(t))=F_0\left(1-\frac tc\right).
\]
The comparison \eqref{eq:higher-comparison}
implies $s\le\bar s(t)$ on $\Sigma_t$.
Using the monotonicity of $j$ and
\eqref{eq:level-flux}, we obtain
\begin{equation}\label{eq:energy-lower-bound}
J
\ge
Q\int_0^c j(\bar s(t))\,dt
=
\frac{cQ}{F_0}I(s_0).
\end{equation}

The elementary identity
\begin{equation}\label{eq:one-dimensional-identity}
nI(s_0)-F_0
=
(n-1)
\left(\sqrt{1+s_0^2}-1\right)s_0^{-1-\beta}
\end{equation}
follows by differentiating both sides with respect
to $s_0$ and observing that both sides tend to
zero as $s_0\downarrow0$.
Combining \eqref{eq:pohozaev-exterior},
\eqref{eq:energy-lower-bound}, and
\eqref{eq:one-dimensional-identity}, we find
\[
\begin{aligned}
nV\left(\sqrt{1+s_0^2}-1\right)
&=nJ-cQ
\ge
\frac{cQ}{F_0}\bigl(nI(s_0)-F_0\bigr)\\
&=
\frac{(n-1)cQ}{F_0}
\left(\sqrt{1+s_0^2}-1\right)s_0^{-1-\beta}.
\end{aligned}
\]
Since $Q=s_0A$, cancellation gives
\begin{equation}\label{eq:reverse-volume-bound}
nV\ge\frac{(n-1)A}{h_0}.
\end{equation}

The curvature bound
\eqref{eq:mean-curvature-lower-bound} makes
$\Gamma$ strictly mean convex.
The Euclidean Heintze--Karcher inequality\cite{HeintzeKarcher1978}
states that
\[
nV\le(n-1)\int_\Gamma\frac1H\,d\sigma.
\]
Together with \eqref{eq:mean-curvature-lower-bound}
and \eqref{eq:reverse-volume-bound}, this gives
\[
nV
\le
(n-1)\int_\Gamma\frac1H\,d\sigma
\le
\frac{(n-1)A}{h_0}
\le
nV.
\]
All inequalities above are therefore equalities,
and $K=B_R(x_0)$ is a ball.

Finally we show that $u$ is radial about $x_0$. 
For an orthogonal transformation $O$, set
\[
u_O(x)=u\bigl(x_0+O(x-x_0)\bigr).
\]
The functions $u_O$ and $u$ solve the same equation
and have the same zero boundary data.
By Lemma~\ref{lem:horizontal-end},
\[
u_O-u\longrightarrow0
\qquad\text{as }|x|\to\infty.
\]
In dimension two, this follows from the common
logarithmic coefficient and constant term;
in higher dimensions, both functions tend to $c$.
The comparison principle on bounded exhaustions
therefore gives $u_O=u$.
Thus $u(x)=U(|x-x_0|)$.

The radial equation is
\[
\left(
r^{n-1}\frac{U'(r)}{\sqrt{1-U'(r)^2}}
\right)'=0.
\]
Since $U'(R)=\kappa$, its first integral is
\[
r^{n-1}\frac{U'(r)}{\sqrt{1-U'(r)^2}}
=
D,
\qquad
D=\frac{\kappa R^{n-1}}{\sqrt{1-\kappa^2}}.
\]
Solving for $U'$ and using $U(R)=0$, we obtain
\[
U(r)
=
\int_R^r
\frac{D}{\sqrt{\rho^{2n-2}+D^2}}\,d\rho.
\]
This is the asserted formula and completes
the proof.
\end{proof}

\section{Rigidity for noncompact connected boundaries}
\label{sec:half-plane-rigidity}

In this section we prove
Theorem~\ref{thm:half-plane-rigidity}. In the bounded-obstacle case considered above, the gradient
vanishes at infinity; although such decay is unavailable in
the present setting, the following lemma shows that the uniform
spacelikeness still follows from the boundary conditions.

\begin{lemma} \label{lem:automatic-uniform-spacelikeness}
Let $\Gamma\subset\mathbb R^2$ be a
$C^{3,\alpha}$ properly embedded curve
diffeomorphic to $\mathbb R$, and let $\Omega$
be one of the two components of
$\mathbb R^2\setminus\Gamma$. Let $\kappa\in(0,1)$ be a constant and we assume that $u$ satisfies \eqref{eq:intro-overdetermined}, then
\[
|Du|\leq\kappa
\quad\text{throughout }\Omega.
\]
\end{lemma}

\begin{proof}
Put
\[
W=\sqrt{1-|Du|^2},
\qquad
P=\frac{1-W}{1+W}
=\frac{|Du|^2}{(1+W)^2}.
\]
Notice that $0\leq P<1$.
On the set $\{Du\neq0\}$, put
$
s=|\nabla_gu|_g=\frac{|Du|}{W}.
$
The two-dimensional identity in
Lemma~\ref{lem:gradient-transform} gives
$
\Delta_gF_2(s)=0,
$
where
$F_2(s)=\log\frac{s}{1+\sqrt{1+s^2}}$. Thus
$P=e^{2F_2(s)}$ and 
\[
\Delta_gP
=4P\,|\nabla_gF_2(s)|_g^2
\geq0
\quad\text{on }\{Du\neq0\}.
\]
The function $P$ is $C^2$ throughout $\Omega$.
At every critical point of $u$, $P$ attains its minimum
value zero, so $\Delta_gP\geq0$ there as well.
Therefore $P$ is subharmonic on all of $\Omega$.

Note 
\[
P=P_\kappa
=\frac{1-\sqrt{1-\kappa^2}}
       {1+\sqrt{1-\kappa^2}} \qquad \text{on } \partial \Omega.
\]
We shall show that $P(x)\leq P_k(x)$ throughout $\Omega$.

To this end, we need a superharmonic barrier as follows. 
For every $x\in\Omega$, choose a nearest boundary
point $y\in\partial\Omega$. Since the segment $[x,y)$ lies
in $\Omega$, the conditions $u(y)=0$ and $|Du|<1$ give
\[
|u(x)|< |x-y|
=\operatorname{dist}(x,\partial\Omega).
\]
Choose $x_*\in\mathbb R^2\setminus\overline{\Omega}$
and $\delta>0$ such that
$
\overline{B_\delta(x_*)}
\subset\mathbb R^2\setminus\overline{\Omega},
$
and write $r(x)=|x-x_*|$. Then
$
\operatorname{dist}(x,\partial\Omega)\leq r(x)-\delta,
$
and consequently
\begin{equation}
\label{eq:lorentzian-radius-lower-bound}
\rho(x):=r(x)^2-u(x)^2
\geq 2\delta r(x)-\delta^2
\geq\delta^2.
\end{equation}

Consider the translated graph immersion
$
X(x)=(x-x_*,u(x))
$ into Lorentz--Minkowski space, so that
$\rho=\langle X,X\rangle_L$.
Let $N$ be its future-directed unit timelike normal.
Maximality gives $\Delta_gX=0$, and hence
\[
\Delta_g\rho=4,
\qquad
\nabla_g\rho=2X^\top.
\]
Since
\[
X^\top=X+\langle X,N\rangle_LN,
\qquad
\langle N,N\rangle_L=-1,
\]
we obtain
\[
|\nabla_g\rho|_g^2
=4\bigl(\rho+\langle X,N\rangle_L^2\bigr).
\]
It follows that
\[
\Delta_g\log\rho
=\frac{\Delta_g\rho}{\rho}
-\frac{|\nabla_g\rho|_g^2}{\rho^2}
=-\frac{4\langle X,N\rangle_L^2}{\rho^2}
\leq0.
\]
Thus the function
$\varphi=\log\frac{\rho}{\delta^2}
$
is nonnegative and superharmonic. Moreover,
\eqref{eq:lorentzian-radius-lower-bound} implies that
$\varphi(x)\to+\infty$ as $|x|\to\infty$ in $\Omega$.

For $\varepsilon>0$, define
$
V_\varepsilon=P-P_\kappa-\varepsilon \varphi.
$
It follows that 
\[
\Delta_gV_\varepsilon\geq0
\quad\text{in }\Omega,
\qquad
V_\varepsilon\leq0
\quad\text{on }\partial\Omega.
\]
The maximum principle, applied on sufficiently large
bounded truncations of $\Omega$, gives
\[
P-P_\kappa\leq\varepsilon \varphi
\quad\text{in }\Omega.
\]

Letting $\varepsilon\downarrow0$ at each fixed point,
we conclude that $P\leq P_\kappa$ in $\Omega$.
Since the function
\[
t\longmapsto
\frac{1-\sqrt{1-t^2}}{1+\sqrt{1-t^2}}
\]
is strictly increasing on $[0,1)$, this implies $|Du|\leq\kappa$. \end{proof}

This uniform spackelikeness enables one to build up a global
conformal parametrization that respects the
noncompact boundary.

\begin{lemma}
\label{lem:global-conformal-parametrization}
Under the assumptions of Theorem~\ref{thm:half-plane-rigidity},
there exists a conformal diffeomorphism
\[
X:\mathbb H\longrightarrow(\Omega,g),
\qquad
\mathbb H=\{z=x+iy:y>0\},
\]
where $g=\delta-du\otimes du$, which extends to a
homeomorphism
$
\mathbb H\cup\mathbb R\cup\{\infty\}
\longrightarrow
\Omega\cup\Gamma\cup\{\infty\}
$
in the spherical topology, with
$
X(\mathbb R)=\Gamma,
$ and $X(\infty)=\infty.$
At every finite boundary point, $X$ extends locally
as a $C^1$ diffeomorphism up to the boundary.
\end{lemma}

\begin{proof}
Since $\Gamma$ is a properly embedded curve
diffeomorphic to $\mathbb R$, the component of
$\mathbb R^2\setminus\Gamma$ other than $\Omega$
contains a closed ball. The boundary conditions give
$Du=\kappa\nu$ on $\Gamma$. Thus
Lemma~\ref{lem:automatic-uniform-spacelikeness} yields
\[
|Du|\leq\kappa<1
\quad\text{in }\Omega,
\]
and consequently
$
(1-\kappa^2)\delta\leq g\leq\delta.
$

The conformal structure determined by $g$ has a
Beltrami coefficient $\mu$ satisfying
\[
|\mu|
\leq
\frac{1-\sqrt{1-\kappa^2}}
     {1+\sqrt{1-\kappa^2}}
<1.
\]
Extend $\mu$ measurably to $\mathbb C$ by setting it
equal to zero outside $\Omega$. By the measurable
Riemann mapping theorem
\cite{Ahlfors2006}, there exists a quasiconformal
homeomorphism
\[
\Phi:\widehat{\mathbb C}\longrightarrow
\widehat{\mathbb C},
\qquad
\Phi(\infty)=\infty,
\]
whose restriction to $\Omega$ is an isothermal
coordinate map for $g$.

Properness of $\Gamma$ implies that
$
\widehat{\Gamma}=\Gamma\cup\{\infty\}
$
is a Jordan curve in $\widehat{\mathbb C}$.
Therefore $\Phi(\widehat{\Gamma})$ is also a Jordan
curve, and
$
\mathcal D=\Phi(\Omega)
$
is one of its complementary domains. The Riemann
mapping theorem and the Carath\'eodory boundary
extension theorem \cite{Pommerenke1992} provide a
conformal map
\[
\psi:\mathbb H\longrightarrow\mathcal D
\]
which extends homeomorphically to the closures in
$\widehat{\mathbb C}$. After precomposing with an
automorphism of $\mathbb H$, we may arrange that
$\psi(\infty)=\infty$. Hence
$
X=\Phi^{-1}\circ\psi
$
has the asserted conformal and topological properties.

It remains to verify regularity at finite boundary
points. This is a local issue. Fix $p\in\Gamma$.
Since $u\in C^2(\overline{\Omega})$ and
$|Du|=\kappa<1$ on $\Gamma$, the metric $g$ admits
a positive definite $C^1$ extension to a neighborhood
of $p$. Local isothermal coordinates for this
extension may be chosen of class $C^{1,\beta}$
for some $\beta\in(0,1)$. In these coordinates,
$\Gamma$ becomes a $C^{1,\beta}$ boundary arc,
and $X$ becomes an ordinary conformal map.

The local boundary regularity theorem of
Kellogg--Warschawski \cite{Pommerenke1992} then
implies that this conformal map extends $C^1$
to the boundary arc, with nonvanishing derivative.
Transforming back to the original coordinates
shows that $X$ extends locally as a $C^1$
diffeomorphism up to $\Gamma$.
\end{proof}

\medskip

We also record a standard consequence of
the representation of positive harmonic
functions on the upper half-plane.

\begin{lemma}\label{lem:positive-harmonic-half-plane}
Let $h$ be a nonnegative harmonic function
on $\mathbb H$, continuous on
$\mathbb H\cup\mathbb R$.
If
\[
h(x,0)=b
\qquad\text{for every }x\in\mathbb R,
\]
where $b\ge0$ is constant, then
$
h(x,y)=b+ay
$
for some $a\ge0$.
\end{lemma}

\begin{proof}
The conclusion is immediate if $h\equiv0$.
Otherwise, the strong maximum principle gives $h>0$
in $\mathbb H$. By the representation theorem for
positive harmonic functions on the upper half-space
\cite[Theorem~7.26]{AxlerBourdonRamey2001}, there exist
$a\geq0$ and a positive Borel measure $\mu$ on $\mathbb R$
such that
\[
h(x,y)
=
ay+\frac1\pi\int_{\mathbb R}
\frac{y}{(x-t)^2+y^2}\,d\mu(t),
\qquad
\int_{\mathbb R}\frac{d\mu(t)}{1+t^2}<\infty.
\]

The approximate identity property of the Poisson kernel,
together with this integrability condition, implies that
\[
h(x,y)\,dx \rightharpoonup \mu
\quad\text{locally as }y\downarrow0.
\]
More precisely, for every $\varphi\in C_c(\mathbb R)$,
\[
\lim_{y\downarrow0}
\int_{\mathbb R}\varphi(x)h(x,y)\,dx
=
\int_{\mathbb R}\varphi(t)\,d\mu(t).
\]
On the other hand, continuity of $h$ on
$\mathbb H\cup\mathbb R$ and the boundary condition
$h(x,0)=b$ imply that $h(x,y)\to b$ uniformly on
each compact interval as $y\downarrow0$. Therefore
\[
\int_{\mathbb R}\varphi(t)\,d\mu(t)
=
b\int_{\mathbb R}\varphi(t)\,dt
\]
for every $\varphi\in C_c(\mathbb R)$, and hence
$\mu=b\,dt$.

Substituting this into the representation formula and
using
\[
\frac1\pi\int_{\mathbb R}
\frac{y}{(x-t)^2+y^2}\,dt=1,
\]
we obtain $h(x,y)=b+ay$, with $a\geq0$.
\end{proof}

Now we are ready to complete the proof. 
\begin{proof}[Proof of
Theorem~\ref{thm:half-plane-rigidity}]
Let $X$ be the conformal parametrization
given by
Lemma~\ref{lem:global-conformal-parametrization},
and write
\[
X^*g=e^{2\varphi}(dx^2+dy^2).
\]
Since $u$ is $g$-harmonic, 
the function
$
U=u\circ X
$
is harmonic on $\mathbb H$.
It is positive in $\mathbb H$ and
continuous with zero boundary values on
the entire real axis.
Lemma~\ref{lem:positive-harmonic-half-plane}
therefore gives
\begin{equation}\label{eq:linear-conformal-height}
U(x,y)=ay
\qquad\text{for some }a>0.
\end{equation}
In particular,
\[
s\circ X
=
|\nabla_{X^*g}U|
=
ae^{-\varphi}>0.
\]
Thus $u$ has no critical points.

Lemma~\ref{lem:gradient-transform} now
applies without truncation and yields
\[
\Delta_g F_2(s)=0
\qquad\text{throughout }\Omega.
\]
Since $F_2(s)<0$ for $s>0$, the function
$
V=-F_2(s)\circ X
$
is positive and harmonic on $\mathbb H$.
The boundary conditions imply
\[
Du=\kappa\nu,
\qquad
s=s_0
\qquad\text{on }\Gamma.
\]
Hence $V$ extends continuously to the
real axis with the constant boundary value
$-F_2(s_0)>0$.
A second application of
Lemma~\ref{lem:positive-harmonic-half-plane}
gives
\[
V(x,y)=-F_2(s_0)+by
\qquad\text{for some }b\ge0.
\]
Combining this identity with
\eqref{eq:linear-conformal-height}, we obtain
$
F_2(s)\circ X
=
F_2(s_0)-\frac ba\,u\circ X.
$
Since $X$ is onto, 
\begin{align} \label{eq:height-gradient-relation}
    F_2(s)=F_2(s_0)-\frac{b}{a} u   \qquad \text{in} \quad \Omega.
\end{align}

Let $H$ denote the signed Euclidean
curvature of $\Gamma$, with the convention
$
H=\operatorname{div}_\Gamma\nu.
$
Thus a circle has positive curvature when
$\nu$ points toward its exterior.
The inward unit conormal in the graph
metric is
$
\eta=\frac{\nu}{\sqrt{1-\kappa^2}},
$
On $\Gamma$ we have
$
\partial_\eta u=s_0.
$
The boundary identity
\eqref{eq:normal-transform}, specialized
to dimension two, gives
$
\partial_\eta F_2(s)=-H.
$
Taking the inward derivative of
\eqref{eq:height-gradient-relation}
therefore yields
\begin{equation}\label{eq:constant-boundary-curvature}
H=\frac{b}{a} s_0
\qquad\text{on }\Gamma.
\end{equation}
In particular, $\Gamma$ has constant
curvature.

A properly embedded curve in $\mathbb R^2$
is complete in its induced Euclidean
metric. A complete connected plane curve
with constant nonzero curvature is a
circle, whereas a complete connected
curve with zero curvature is a line.
Since $\Gamma$ is noncompact and
diffeomorphic to $\mathbb R$, the circle
case is excluded.
Thus
$
H=0.
$
Because $s_0>0$, equation
\eqref{eq:constant-boundary-curvature}
implies $b=0$.
Consequently,
\[
F_2(s)=F_2(s_0)
\qquad\text{in }\Omega.
\]
The function $F_2$ is strictly increasing,
so
$
s\equiv s_0,$ and $
|Du|\equiv\kappa.$

We finally show that $u$ is affine.
Differentiating $|Du|^2=\kappa^2$ gives
$
D^2u\,Du=0.
$
The nondivergence form of the maximal
surface equation is
\[
\Delta u+
\frac{D^2u(Du,Du)}{1-|Du|^2}=0,
\]
and therefore $\Delta u=0$.
The Euclidean Bochner identity gives
\[
0=\frac12\Delta|Du|^2
=
|D^2u|^2+Du\cdot D(\Delta u)
=
|D^2u|^2.
\]
Hence $D^2u=0$ and $u$ is affine.

We already know that $\Gamma$ is a line,
and $\Omega$ is one of its complementary
half-planes.
After a Euclidean rigid motion we may
write
\[
\Gamma=\{x_2=0\},
\qquad
\Omega=\{x_2>0\}.
\]
The boundary conditions $u=0$ and
$\partial_\nu u=\kappa$ then identify
the affine function as
$
u(x_1,x_2)=\kappa x_2.
$
\end{proof}

\end{document}